\documentclass[leqno]{amsart}
\usepackage{amsmath,amsfonts,amsthm,amssymb, mathrsfs, mathabx, accents, mathdots, textcomp, fontenc,  vntex}
\usepackage[bookmarks]{hyperref}
\hypersetup{backref, colorlinks=true}
\usepackage[french, english]{babel}
\usepackage{amsmath,amsfonts,amsthm,amssymb, mathrsfs, mathabx, accents, mathdots}
\usepackage[all]{xy}
 \usepackage{graphicx}
\usepackage[all]{xy}
\newtheorem{theorem}{Th\'eor\`eme}[section]
\newtheorem{lemma}[theorem]{Lemme}
\newtheorem{proposition}[theorem]{Proposition}

\newtheorem{remark}[theorem]{Remarque}

\newtheorem{example}[theorem]{Exemple}
\newtheorem{definition}[theorem]{D\'efinition}

\newenvironment{equationth}{\stepcounter{theorem}\begin{equation}}{\end{equation}}
\newenvironment{preuve}{{\em{\noindent \textbf{Preuve.} }}}
{\hfill $\blacksquare$}

\def\C{ \mathbb{C}}

\def\R{ \mathbb{R}}

\def\N{ \mathbb{N}}

\def\rond{\mathaccent"7017}

\begin{document}

\large 

\title[]{Une stratification de Thom-Mather de l'ensemble asymptotique  d'une application polynomiale $F: \C^n \to \C^n$} 
\makeatother

\author[Nguy\~{\^e}n Th\d{i} B\'ich Th\h{u}y]{Nguy\~{\^e}n Th\d{i} B\'ich Th\h{u}y}
\address[{Nguy\~{\^e}n Th\d{i} B\'ich Th\h{u}y}]{UNESP, Universidade Estadual Paulista, ``J\'ulio de Mesquita Filho'', S\~ao Jos\'e do Rio Preto, Brasil}
\email{bichthuy@ibilce.unesp.br}
\maketitle \thispagestyle{empty}
\begin{abstract}
Let $F: \C^n \to \C^n$ be a polynomial mapping. We prove that the stratification of the asymptotic set of $F$  defined by {\it ``la m\'ethode des fa{\c c}on''} in \cite{Thuy} is a Thom-Mather stratification. 

\medskip 

\medskip 

\noindent R\'ESUM\'E. Soit  $F: \C^n \to \C^n$ une application polynomiale. Nous prouvons que la stratification de l'ensemble asymptotique  de $F$ d\'efinie par {\it ``la m\'ethode des fa{\c c}ons''} in \cite{Thuy} est une stratification de Thom-Mather. 
\end{abstract}

\begin{center} 
{\large \bf INTRODUCTION}
\end{center}

\medskip

Soit $F : \C^n_{(x)} \to \C^n_{(\alpha)}$ une application polynomiale. L'ensemble asymptotique de $F$, not\'e $S_F$, est l'ensemble des points du but en lesquels l'application $F$ n'est pas propre. 
Dire que l'application $F$ n'est pas propre en un point  du but $a \in  \C^n_{(\alpha)}$ peut \^etre caract\'eris\'e des deux mani\`eres suivantes~: 

1) Il existe une suite  $\{ \xi_k\}_{k \in \N} \subset \C^n_{(x)}$ dans la source telle que $\{ \xi_k\}$ tende vers l'infini et telle que l'image $F(\xi_k)$ tende vers $a$. 

\noindent Ici, ``la suite $\{ \xi_k\}$ tend vers l'infini'' signifie que la norme euclidienne $\vert \xi_k \vert$ de $\xi_k$ dans $\C^n_{(x)}$ tend vers l'infini. 

2) Il existe une courbe diff\'erentiable  
$$\gamma : (0, +\infty) \, \to \C^n_{(x)}, \quad \gamma(u) = (\gamma_1(u), \ldots, \gamma_n(u))$$ 
tendant vers l'infini et telle que $F \circ \gamma (u)$ tend vers $a$ lorsque $u$ tend vers l'infini.

Dans les ann\'ees 90, Jelonek a \'etudi\'e l'ensemble asymptotique associ\'e \`a une application polynomiale $F: \C^n \to \C^n$ de mani\`ere approfondie et il en a d\'ecrit les  principales propri\'et\'es \cite{Jelonek1}. 
La compr\'ehension de la structure de cet ensemble est tr\`es importante par sa relation avec la  Conjecture  Jacobienne (voir, par exemple, \cite{Essen}).  


Dans \cite{Thuy}, nous fournissons une m\'ethode, appell\'e {\it la m\'ethode des fa{\c c}ons} pour stratifier l'ensemble asymptotique $S_F$. La stratification obtenue est une stratification diff\'erentiable et satisfait la condition de fronti\`ere.
 Cette m\'ethode \'etudie les comportaments 
des courbes tendant vers l'infini et telles que leurs images tendent vers les points de $S_F$ comme suit~: ces suites seront  labellis\'ees sous la forme de ``fa{\c c}ons (\'etoile)'' (D\'efinition 2.13 de \cite{Thuy}) telles que chaque fa{\c c}on ``\'etoile'' diff\'erente d\'efinit des  courbes co\-rres\-pon\-dantes parall\`eles localement, qui, en fait, d\'efinissent un feuilletage complexe de dimension 1 (Proposition 2.18 de \cite{Thuy}) dans l'espace source $\C^n$. Les images des feuilletages  co\-rres\-pon\-dant aux  diff\'erentes fa{\c c}ons d'approcher la partie  singuli\`ere de $S_F$, ce qui nous permet de d\'ecomposer celle-ci en strates. Nous obtenons une stratification diff\'erentiable de l'ensemble asymptotic $S_F$, qui satisfait aussi la condition de fronti\`ere (Th\'eor\`eme 4.1 de \cite{Thuy}).

  Nous montrons dans cet article que la stratification d\'efinie par les {\it fa{\c c}ons} dans \cite{Thuy} est une stratification de Thom-Mather. 
Pour montrer cela, nous utilisons l'autre caract\'erisation 
de l'ensemble asymptotique, comme limite des i\-ma\-ges de courbes, 
lesquelles tendent vers l'infini. 
Lorsque ces courbes sont parall\`eles (d\'efinissent un feuilletage), la limite des i\-ma\-ges est une strate, 
lisse, obtenue comme limite des images du feuilletage transverse. Nous monstrons que les images de telles 
courbes constituent les rayons des voisinages tubulaires  ``\`a la Thom-Mather'' ce qui d\'emontre notre Th\'eor\`eme \ref{theoremThom-Mather}. 

Une autre motivation de cette \'etude r\'eside dans le fait que conna\^itre une stratification d'une vari\'et\'e singuli\`ere permet de calculer son homologie d'intersection: Dans \cite{Valette}, les  auteurs ont contruit des vari\'et\'es singuli\`eres $V_F$ associ\'ees \`a une application polynomial $F: \C^n \to \C^n$ telle que l'homologie d'intersection des ces vari\'et\'es $V_F$  caract\'erise  la propret\'e de $F$ dans le cas le jacobian de $F$ est partout non nul (Th\'eor\`eme 3.2 de \cite{Valette} et Th\'eor\`eme 4.5 de \cite{Thuy2}). De plus, une stratification de l'ensemble asymptotique $S_F$ permet de fournir une stratification de la vari\'et\'e $V_F$. Notons que pour une application fix\'ee $F$, nous avons beaucoup de vari\'et\'es associ\'ees $V_F$, mais les r\'esultats de l'homologie d'intersections des vari\'et\'es $V_F$ dans \cite{Valette} et dans \cite{Thuy2} ne changent pas. 
Nous savons aussi qu'en general, l'homologie d'intersection d'une vari\'et\'e singuli\`ere d\'epend de sa stratification. Cependant, l'homologie d'intersection d'une vari\'et\'e singuli\`ere est invariant avec une stratification de Thom-Mather. Alors, une stratification de Thom-Mather de l'ensemble asymptotique $S_F$ est n\'ecessaire pour calculer l'homologie d'intersection des vari\'et\'es $V_F$. Nous pr\'ecisons donc les r\'esultats obtenus dans \cite{Valette} et dans \cite{Thuy2}.

\medskip

\begin{center}
{ \large \bf NOTATION}
\end{center}

\medskip

Nous consid\'erons dans cet article des applications polynomiales $F: \C^n \to \C^n$. Nous \'ecrivons souvent $F: \C^n_{(x)} \to \C^n_{(\alpha)}$ pour distinguer entre la source et le but. 

\medskip

\medskip

\begin{center}
{ \large \bf 0. PR\'ELIMINAIRES}
\end{center}

\medskip

\subsection{Stratifications} \label{sectionStratification}

\begin{definition}
{\rm 
 {Soit $V$ une vari\'et\'e (diff\'erentiable ou alg\'ebrique, ou analytique) de dimension $m$. {\it Une stratification}  (${{\mathscr{S}}}$) de $V$ est la donn\'ee d'une filtration 
$$V=V_m \supseteq V_{m-1} \supseteq V_{m-2} \supseteq \dots \supseteq V_1 \supseteq V_0 \supseteq V_{-1} = \emptyset $$
de $V$ telle que toutes les diff\'erences $X_{i} = V_{i} \setminus V_{i-1}$ sont ou bien   vides ou bien unions localement finies de sous-vari\'et\'es  lisses connexes et localement ferm\'ees de dimension $i$, appel\'ees strates.

%

Soit $S_i$ une strate de $V$ et $\overline{S_i}$ est son adh\'erence dans $V$. 
Si $\overline{S_i} \setminus S_i$ est l'union de strates de $V$, pour toute strate $S_i$ de $V$, 
alors nous disons que la stratification de $V$ satisfait la propri\'et\'e de fronti\`ere. 
}
}
\end{definition}

\begin{definition} [voir \cite{Thom}, \cite{Mather1}] \label{definitionthommather}
{\rm Soit $V$ une sous-vari\'et\'e d'une vari\'et\'e lisse $M$. Nous disons qu'une stratification de $V$ est {\it une stratification de Thom-Mather} si chaque strate $S_i$ est une vari\'et\'e diff\'erentiable de classe $\mathcal C^{\infty}$  et si pour chaque strate $S_i$ nous avons

a) $ \,$ un voisinage ouvert (voisinage tubulaire) $T_i$ de $S_i$ dans $M$,

b) $ \,$ une r\'etraction continue $\pi_i$ de $T_i$ sur $S_i$,

c) $ \,$ une fonction continue (tubulaire) $\rho_i: T_i \to [0, +\infty)$ qui est $\mathcal C^\infty$ sur la partie r\'eguli\`ere de $V \cap T_i$,

\noindent tels que $S_i = \{ x \in T_i  : \rho(x) = 0 \}$ et si $S_i \subset \overline {S_j}$, alors

i) $ \,$ l'application restriction $(\pi_i, \rho_i): T_i \cap S_j \to S_i \times [0, +\infty)$ est une immersion lisse,

ii) $ \,$ pour $x \in T_i \cap T_j$ tel que $\pi_j(x) \in T_i$, nous avons les relations de commutation :

$\qquad$ 1) $\pi_i \circ \pi_j(x) = \pi_i(x),$ 

$\qquad$ 2) $\rho_i \circ \pi_j(x) = \rho_i(x),$ 

\noindent lorsque les deux membres de ces formules sont d\'efinis.
}
 \end{definition}

%
%
%
%
%

\begin{remark}
{\rm Les conditions de Whitney impliquent des conditions de Thom-Mather. }
 \end{remark}

\subsection{L'ensemble asymptotique.} \label{ensembleJelonek}

Soit $F : \C^n_{(x)} \to \C^n_{(\alpha)}$ une application polynomiale. Notons $S_F$ l'ensemble des points du but pour lesquels l'application $F$ n'est pas propre, {\it i.e.},  
$$S_F = \{ a \in \C^n_{(\alpha)} \text{ tel que } \exists \{ \xi_k\}_{k \in \N} \subset \C^n_{(x)}, \vert \xi_k \vert  \text{ tend vers l'infini et } F(\xi_k) \text{ tend vers } a\},$$
o\`u $ \vert \xi_k \vert$ est  la norme euclidienne de $\xi_k$ dans $\C^n$. 
L'ensemble $S_F$ est appel\'e l'ensemble asymptotique de $F$.

\begin{lemma} \label{lemmecourbe} 
{\rm Le point $a$ appartient \`a $S_F$ si et seulement s'il existe une courbe diff\'erentiable  $\gamma(u) : (0, +\infty) \, \to \C^n_{(x)}$ tendant vers l'infini et telle que $F \circ \gamma (u)$ tend vers $a$ lorsque $u$ tend vers l'infini. 
}
\end{lemma}

Rappelons qu'il suffit, pour d\'efinir $S_F$ de consid\'erer des courbes $\gamma(u) = (\gamma_1(u), \ldots, \gamma_n(u))$ tendant vers l'infini, au sens suivant ~: chaque coordonn\'ee $\gamma_1(u), \ldots, \gamma_n(u)$ de cette courbe ou bien tend vers l'infini ou bien converge. C'est ce que nous ferons dans cet article.

\begin{definition} \label{definitiondominant}
{\rm Une application polynomiale $F : \C^n_{(x)} \to \C^n_{(\alpha)}$ est dite {\it dominante} si l'adh\'erence de $F(\C^n_{(x)})$ est dense dans $\C^n_{(\alpha)}$, 
c'est-\`a-dire $\overline{F(\C^n_{(x)})} = \C^n_{(\alpha)}$. 
}
\end{definition}


\begin{theorem} \cite{Jelonek1} \label{theoremjelonek1}
Soit $F = (F_1, \ldots, F_n) : \C^n \rightarrow \C^n$ une application polynomiale dominante. Pour chaque $ i = 1, \ldots, n$, consid\'erons une \'equation irr\'eductible de $x_i$ sur $\C[F_1, \ldots , F_n]$ :
$$\sum_{k=0}^{n_i}\phi_k^i(F)x_i^{n_i-k}=0,$$
o\`u les $\phi_k^i$ sont des polyn\^omes et les $n_i$ sont des entiers positifs. Alors nous avons 
$$S_F = \bigcup_{i=1}^n \{a \in \C^n: \phi_0^i(a)=0 \}.$$

\end{theorem}

\section{La m\'ethode  des fa{\c c}ons}

Dans \cite{Thuy}, nous fournissons  {\it la m\'ethode des fa{\c c}ons} pour stra\-ti\-fi\-er la vari\'et\'e asymptotique d'une application polynomiale dominante $F: \C^n \to \C^n$. 
 La d\'etermination des strates de la stratification de l'ensemble asymptotique n\'ecessite deux  \'etapes. 
La premi\`ere \'etape, qui, bien que est tr\`es significative pour la construction de la m\'ethode,  fournit une partition dont les \'el\'ements peuvent \^etre des 
vari\'et\'es singuli\`eres. Un raffinement de cette partition est n\'ecessaire, en utilisant une version 
plus fine: les fa\c cons  ``\'etoile''. 
 De mani\`ere plus pr\'ecise, les courbes $\gamma$ tendant vers l'infini et telles que leurs images $F \circ \gamma$ tendent vers les points de $S_F$ seront  labellis\'ees sous la forme de ``{\it fa{\c c}ons}'' et ``{\it fa{\c c}ons \'etoile}'', respectivement. 

\subsection{Fa{\c c}ons}
 Dans la premi\`ere \'etape, nous d\'efinissons ``{\it une fa{\c c}on du point $a \in S_F$}''. 
Un point $a$ de $S_F$ est la limite de $ F\circ \gamma (u)$, o\`u $\gamma(u) = (\gamma_1(u), \ldots, \gamma_n(u)) : (0, +\infty) \, \to \C^n_{(x)}$ est une courbe  tendant vers l'infini. Nous classons les 
coordonn\'ees $\gamma_1(u), \ldots, \gamma_n(u)$ de la courbe $\gamma (u)$ en trois cat\'egories~: 
i) les coordonn\'ees $\gamma_{i_r}(u)$ tendant vers l'infini (cette cat\'egorie n'est pas vide);
ii) les coordonn\'ees $\gamma_{j_s}(u)$ telles que $\lim_{u \to \infty} \gamma_{j_s}(u)$ est un nombre complexe  ``ind\'ependant du 
point $a$ dans un voisinage de $a$ dans $S_F$''. 
Cela signifie qu'il existe des points $a'$ voisins de $a$ dans $S_F$ 
et des courbes $\gamma^{a'} = (\gamma^{a'}_1(u), \ldots, \gamma^{a'}_n(u)): (0, +\infty) \, \to \C^n_{(x)}$ tendant vers l'infini telles que $\lim_{u \to \infty} F \circ \gamma^{a'}(u) = a'$ et $\lim_{u \to \infty} \gamma^{a'}_{j_s} (u) = \lim_{u \to \infty} \gamma_{j_s} (u) = constante$; 
iii) les coordonn\'ees $\gamma_{ i_l} (u)$ telles que $\lim_{u \to \infty} \gamma_{i_l}(u)$ est un nombre complexe  ``d\'ependant du point $a$'' (ce cas est le cas contraire du cas ii)).
 
Nous d\'efinissons ``{\it une fa{\c c}on du point $a \in S_F$}''
comme un $(p,q)$-uple $(i_1, \ldots , i_p)[j_1, \ldots, j_q]$ d'entiers o\`u  les entiers $i_1, \ldots , i_p$ sont
les indices des coordonn\'ees de la premi\`ere cat\'egorie et $j_1, \ldots, j_q$ 
les indices des coordonn\'ees de la seconde cat\'egorie. Pour comprendre mieux la d\'efinition de {\it fa{\c c}ons}, voir la D\'efinition  1.2 et la section 1.2 de \cite{Thuy}. 
Nous donnons ici un exemple pour \'eclairer cette D\'efinition.

 \begin{example} \label{ex1}
{ \rm Soit l'application $F: \C^3_{(x_1,x_2,x_3)} \to \C^3_{(\alpha_1, \alpha_2,\alpha_3)}$ d\'efinie par 
$$F(x_1,x_2,x_3) = (x_1x_2,x_2x_3,x_1x_2x_3).$$ Alors l'ensemble $S_F$ est \'egal \`a ${(S_F)}_1 \cup {(S_F)}_2 \cup {(S_F)}_3$,
 o\`{u} ${(S_F)}_1 = \{ \alpha_2 = 0 \} $, ${(S_F)}_2 = \{ \alpha_3 = 0 \} $
 et ${(S_F)}_3 = \{ \alpha_1 = 0 \}.$ En fait   

i) Pour un point $a = (\alpha_1 , 0, \alpha_3) \in {(S_F)}_1$ tel que $\alpha_1 \neq 0$ et $ \alpha_3 \neq 0$, il existe une courbe $\gamma = \left( \alpha_1 u, \frac{1}{u},  \frac{\alpha_3}{\alpha_1} \right)$  tendant vers l'infini  telle que l'image $F \circ \gamma $ tende vers $a$. 
Nous disons qu'une  {\it ``fa{\c c}on''} du point $a$ est  $(1)[2]$, o\`u
\begin{enumerate}
\item Le symbole ``(1)''  signifie que la premi\`ere coordonn\'ee $\gamma_1(u) = \alpha_1 u$ de la courbe $\gamma(u)$ tend vers l'infini. 
\item Le symbole ``[2]'' signifie que  
la deuxi\`eme coordonn\'ee $\gamma_2(u) = \frac{1}{u}$ de la courbe $\gamma(u)$ tend vers z\'ero, qui est un nombre complexe fix\'e, lequel ne d\'epend pas du point $a = (0, \alpha_2, \alpha_3)$ quand $a$ d\'ecrit  ${(S_F)}_2 = \{ \alpha_3 = 0 \} $. 

\item La troisi\`eme coordonn\'ee $\gamma_3(u) =  \frac{\alpha_3}{\alpha_1} $ de la courbe $\gamma(u)$ tend vers un nombre complexe 
 $\frac{\alpha_3}{\alpha_1}$ d\'ependant du point  $a = (\alpha_1 , 0, \alpha_3)$  quand $a$ varie. 
 Alors, avec notre choix de labeliser les courbes, l'indice ``2'' n'apparaitra pas dans la façon (1)[2]
\end{enumerate}

Notons que  nous pouvons v\'erifier facilement que toute courbe $\hat{\gamma} (u)$ tendant vers l'infini telle que $F \circ \hat{\gamma} (u)$ tende vers un point de ${(S_F)}_1$ admet la m\^eme fa{\c c}on $(1)[2]$. Nous disons que  la fa{\c c}on de ${(S_F)}_1$ est $(1)[2]$. 

De la m\^eme mani\`ere que dans le cas i), nous avons : 

ii) Pour un point $a = (\alpha_1 , \alpha_2 , 0) \in {(S_F)}_2$ 
tel que $\alpha_1 , \alpha_2 \neq 0$, 
il existe une courbe 
$\gamma = \left( \frac{\alpha_1}{u}, u  , \frac{\alpha_2}{u} \right)$ 
telle que $\gamma$ tend vers l'infini et telle que $F \circ \gamma$ tend vers $a$. 
Nous avons la fa{\c c}on de ${(S_F)}_2$ est  $(2)[1,3]$.

iii) Pour un point $a = (0, \alpha_2 , \alpha_3 ) \in {(S_F)}_3$, 
tel que $\alpha_2, \alpha_3 \neq 0$, 
il existe une courbe  $\gamma = \left( \frac{\alpha_3}{\alpha_2}, \frac{\alpha_2}{u} , u \right) $ 
 telle que $\gamma$ tend vers l'infini et telle que 
$F \circ \gamma$ tend vers $a$. La fa{\c c}on de ${(S_F)}_3$ est donc 
$(3)[2]$.

iv) Pour un point $a = ( \alpha_1, 0,0 ) \in {(S_F)}_1 \cap {(S_F)}_2$, o\`u $\alpha_1 \neq 0$, 
il existe deux courbes  $\gamma =  \left( \alpha_1u, \frac{1}{u}, \frac{1}{u} \right)$ 
et  
$\hat{\gamma} = \left( \frac{1}{u}, u \alpha_1 , \frac{1}{u^2}  \right)$ 
telles que $\gamma$ et $\hat{\gamma}$ tendent vers l'infini et  
leurs images $F \circ \gamma$ et $F \circ \hat{\gamma}$ 
tendent vers $a$.
L'ensemble $0 \alpha_1 \setminus \{ 0 \}$ a donc deux fa{\c c}ons  (1)[2,3] et  (2)[1,3].

v) Pour un point $a = (0, \alpha_2, 0 ) \in S_{F_2} \cap S_{F_3} $,  o\`u $\alpha_2 \neq 0$, 
il existe deux courbes  $\gamma = \left(0, u, \frac{\alpha_2}{ u}\right)$ et  
$ \hat{\gamma}= \left( 0, \frac{\alpha_2}{u} ,  u \right)$ telles que telles que $\gamma$ et $\hat{\gamma}$ tendent vers l'infini et  
leurs images $F \circ \gamma$ et $F \circ \hat{\gamma}$ 
tendent vers $a$.
 L'ensemble $0 \alpha_2 \setminus \{ 0 \}$ a donc deux fa{\c c}ons (2)[1,3] et (3)[1,2]. 

vi) Pour un point $a = (0 , 0, \alpha_3 ) \in {(S_F)}_3 \cap {(S_F)}_1$, o\`u $\alpha_3 \neq 0$, 
il existe une courbe $\gamma = \left( u, \frac{1}{ u^2} , \alpha_3  u \right)$ telle que $\gamma$ tend vers l'infini et telle que $F \circ \gamma$ tend vers $a$. 
 La fa{\c c}on de l'ensemble $0 \alpha_3 \setminus \{ 0 \}$ est donc 
$ (1,3)[2].$   

vii) Enfin, consid\'erons l'origine $0 \in {(S_F)}_1 \cap {(S_F)}_2 \cap {(S_F)}_3$, 
il existe quatre courbes 
$\left(  \lambda, \frac{1}{u^2} , u \right)$ avec $\lambda \in \C$, 
 $(u, 0 , u)$, $\left( \frac{1}{u^2} , u, 0 \right)$ et  
$( u, 0 , \mu ) \}]$ avec $\mu \in \C$  tendent vers  l'infini et leurs images tendent vers $a$. 
L'origine a donc quartre fa{\c c}ons (3)[1,2], (1,3)[2], (2)[1,3] et (1)[2,3].

La partition de l'ensemble $S_F$, d\'efinie par des fa{\c c}ons, est donn\'ee par la filtration :
$$ S_F \supset 0\alpha_1 \cup 0\alpha_2 \cup 0\alpha_3  \supset \{ 0 \} \supset \emptyset. $$

\begin{figure}[h]\begin{center}
\includegraphics[scale=0.7]{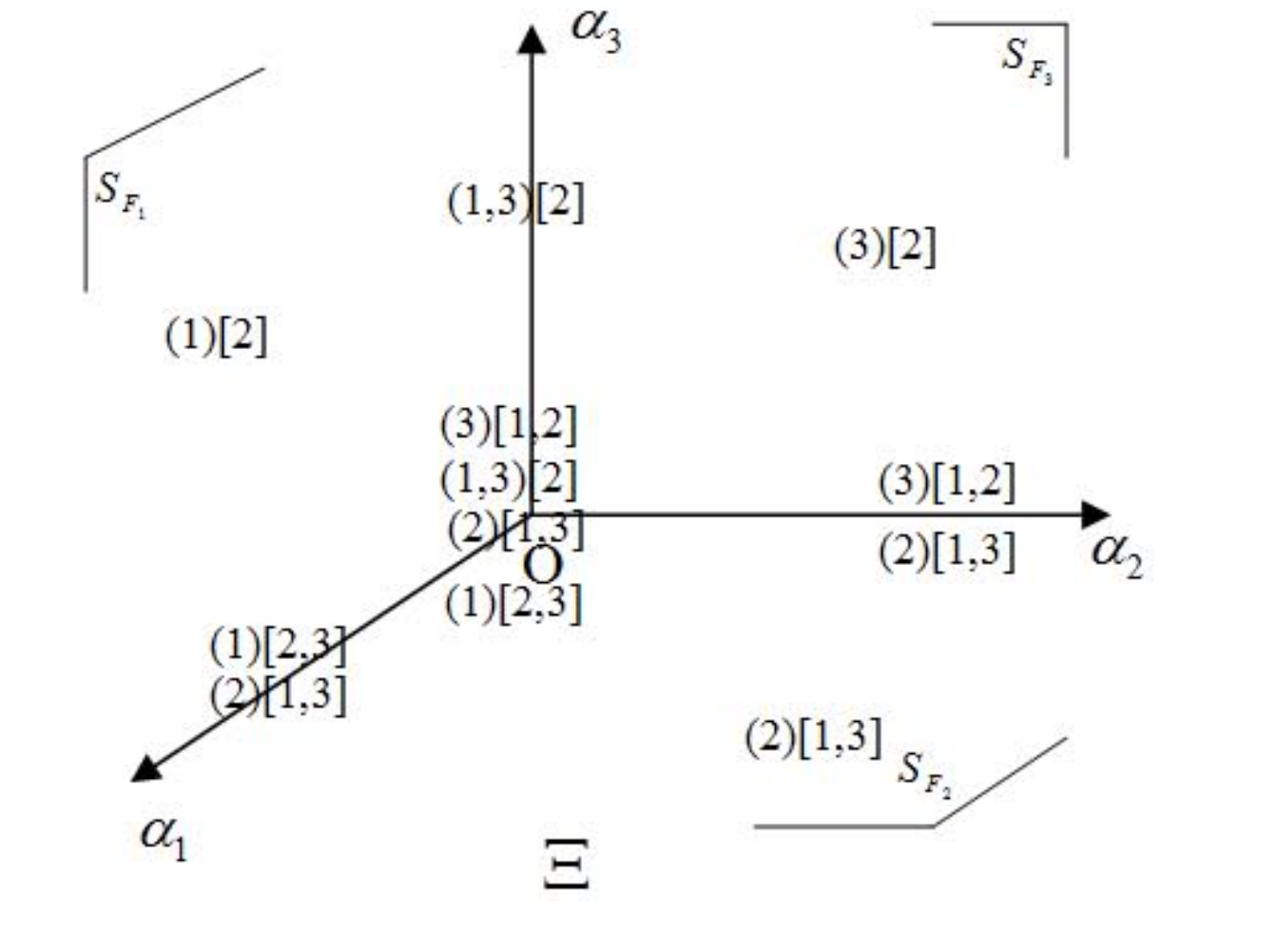}
\caption{La partition de $S_F$ d\'efinie par $\Xi$, pour l'application polynomiale dominante  $F(x_1,x_2,x_3) = (x_1x_2,x_2x_3,x_1x_2x_3)$}
\label{figureex1}
\end{center}\end{figure}
}
\end{example}

\footnote{Revoir l'exemple, s'il est trop compliqu\'e? On a vraiment besoin un exemple comme \c ca pour \'eclairer l'article (oui)? Si oui, refaire le dessin.}

Notons $\Xi(a)$ l'ensemble de toutes les fa{\c c}ons du point $a$. La partition de $S_F$ d\'efinie par la relation 
\begin{equationth} \label{partition facon}
a_1 \sim a_2 \text{ si et seulement si } \Xi(a_1) = \Xi(a_2)
\end{equationth}
est une partition finie de $S_F$ (Propostion 1.6 de \cite{Thuy}).

Les sous-vari\'et\'es de la partition de l'ensemble asymptotique $S_F$ d\'efinie par la relation \ref{partition facon}, bien que tr\`es significatives pour la construction de {\it la m\'ethode des fa{\c c}ons},   se r\'ev\`elent  malheureusement 
pouvoir \^etre singuli\`eres, comme le montre de l'exemple suivant (Proposition 1.11 de \cite{Thuy}). 

\begin{example}  \label{contre-exemple}
 {\rm Soit $F : \C^2_{(x_1, x_2)} \to \C^2_{(\alpha_1, \alpha_2)}$ l'application polynomiale dominante telle que 
$$F(x_1,x_2) = \left({(x_1x_2)}^2, {(x_1x_2)}^3 + x_1 \right).$$ 
Si la courbe $\gamma(u) = (\gamma_1(u), \gamma_2(u)) :  (0, +\infty) \rightarrow  \C^2_{(x_1, x_2)} $ 
tend vers l'infini telle que 
$F \circ \gamma (u) = \left( { (\gamma_{1}(u) \gamma_2 (u))}^2, {(\gamma_{1}(u) \gamma_{2}(u))}^3 + \gamma_{1}(u) \right)$ 
ne tend pas vers l'infini,  alors $\gamma_{1}(u)$ ne peut pas tendre vers l'infini. 
Comme $\gamma(u)$ tend vers l'infini, alors $\gamma_ 2 (u)$ doit tendre vers l'infini et donc $S_F$ n'admet  qu'une seule fa{\c c}on $\kappa = (2)[1]$. Si nous choisissons les courbes coordonn\'ees $\gamma_1(u)$ et $\gamma_2(u)$ tendant vers 0 et l'infini, respectivement, telles que le produit $\gamma_1(u) \gamma_2(u)$ tend vers un nombre complexe $\alpha \in \C$, alors $F\circ \gamma (u)$ tend vers  $(\alpha^2, \alpha^3)$. 
 L'ensemble $S_F$  est donc la courbe $ \alpha_1^3 = \alpha_2^2$ dans $\C^2_{(\alpha_1, \alpha_2)}$ ayant un point singulier \`a l'origine (voir Figure \ref{facon1}). 

\begin{figure}[h!]
\centering
\includegraphics[scale=0.7]{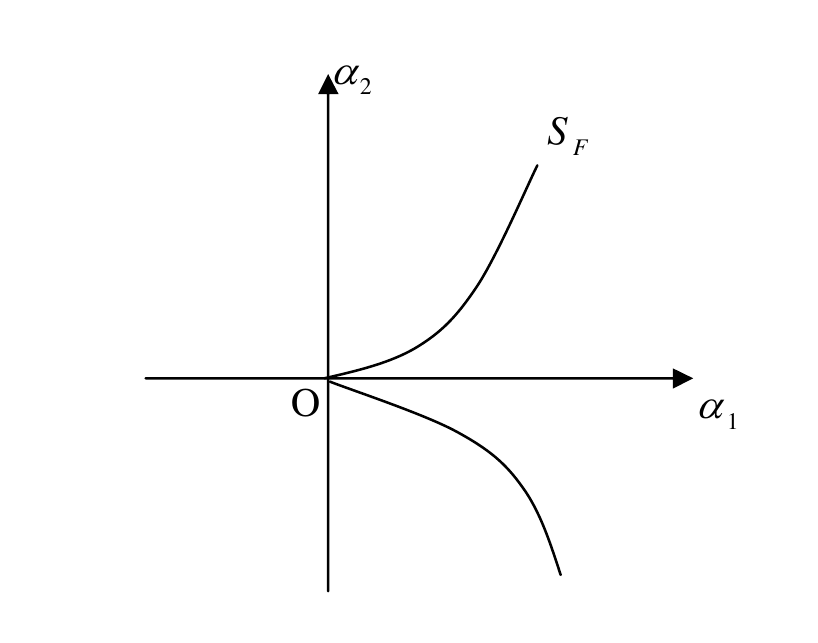}
\caption{L'ensemble asymptotique $S_F$ de l'application $F(x_1,x_2) = \left({(x_1x_2)}^2, {(x_1x_2)}^3 + x_1 \right).$}
\label{facon1}
\end{figure}
}
\end{example}

Donc, un raffinement de la d\'efinition des  fa{\c c}ons  est n\'ecessaire afin d'obtenir une stra\-ti\-fi\-ca\-tion en strates lisses. Cela est la deuxi\`eme \'etape de ``La m\'ethode des fa{\c c}ons''.

\section{Fa{\c c}ons \'etoile}

La question est donc de savoir pourquoi, dans l'exemple de l'exemple \ref{contre-exemple},   
la partition de $S_F$ d\'efinie par les fa{\c c}ons n'est pas  une stratification. Nous devons comprendre la diff\'erence entre ``la fa{\c c}on $[2](1)$ \`a l'origine'' et ``la fa{\c c}on $[2](1)$ au point $a \in S_F \setminus \{ 0 \}$''~:
Pour un point $a = (\alpha^2, \alpha^3)$ dans $S_F \setminus \{ 0 \}$, nous devons choisir la  courbe $ \left( \frac{1}{u^r}, \alpha u^s \right)$, 
o\`u $r$ doit \^etre \'egal \`a $s$.  \`A l'origine 0, nous devons choisir une courbe $ \left( \frac{1}{u^r},  u^s \right)$ o\`u $r > s$.  Nous voyons que nous devons \'etudier la vitesse de  tendre vers zero et  de tendre vers l'infini des courbes coordonn\'ees. Dans l'exemple ci-dessus, si la vitesse de tendre vers zero  de la premi\`ere courbe coordonn\'ee est plus rapide que la vitesse de tendre vers l'infini  de la deuxi\`eme courbe coordonn\'ee, nous tendons \`a la partie singuli\`ere de $S_F$.

Pour  formaliser l'id\'ee des {\it fa{\c c}ons \'etoile}, nous avons donc besoin d'abord de d\'ecrire  ``la vitesse de  tendre vers zero et  de tendre vers l'infini des courbes coordonn\'ees dans $\C$''.

1) D'abord, nous d\'ecrivons ``la vitesse de  tendre vers zero et  de tendre vers l'infini des courbes coordonn\'ees dans $\C$'': Notons que si une courbe $\rho :  (0, +\infty) \, \to \C$
 tend vers un nombre complexe $\lambda$ dans $\C$, alors 
nous pouvons consid\'erer que cette courbe tend vers $0$ par changement de  variables $\rho - \lambda$. Pour cette raison, nous pouvons supposer qu'une courbe $\rho$ dans $\C$ ou bien tend vers zero, ou bien tend vers l'infini. Nous disons que 

+ la courbe $\rho$ tend vers 0 avec le degr\'e $t$ si $\rho(u) = {(\lambda/u)}^t + \cdots$, o\`u $\lambda \in \C \setminus \{ 0 \}$, $t \in \R, t > 0$ et les \'el\'ements dans ``$\ldots$'' sont de la forme ${(\lambda'/u)}^r$ avec $r > t $. 

+ la courbe $\rho$ tend vers l'infini avec le degr\'e $t$ si $\rho (u) =\lambda u^t + \cdots$, 
o\`u  $\lambda \in \C \setminus \{ 0 \}$, 
$t \in \R, t > 0$ et les \'el\'ements dans ``$\ldots$'' sont de la forme 
$\lambda' k^r$ avec $ r < t$. 

2) Comparer la vitesse de tendre vers zero et de tendre vers l'infini des courbes coordonn\'ees $\gamma_1, \ldots, \gamma_n$ de la courbe $\gamma = (\gamma_1, \ldots, \gamma_n) :  (0, +\infty) \, \to \C^n_{(x)}$ comme suit : 
{Soit   $\kappa = (i_1, \ldots, i_p)[j_1, \ldots, j_q]$ une fa{\c c}on de $S_F$ et soit $\gamma = (\gamma_1, \ldots, \gamma_n) :  (0, +\infty) \, \to \C^n_{(x)}$  une courbe tendant vers l'infini avec la 
fa{\c c}on $\kappa$.  Supposons que~: 

$\quad \quad \quad$ $\gamma_{i_r}(u)$ tend vers l'infini avec le degr\'e $l_{i_r}$, pour $r=1, \ldots, p$,

$\quad \quad \quad$ $\gamma_{j_s}(u)$ tend vers 0 avec le degr\'e $l_{j_s}$, pour $s=1, \ldots, q$. 

\noindent Le $(p+q)$-uple $(l_{i_1}, \ldots, l_{i_p}, l_{j_1}, \ldots, l_{j_q}) $ est appel\'e {\it $(p+q)$-uple associ\'e \`a la courbe}  
$\gamma$. 

Notons que,  par un changement de param\`etre, nous povons toujours supposer que les nombres $l_{i_1}, \ldots, l_{i_p}, l_{j_1}, \ldots, l_{j_q} $ du  $(p+q)$-uple  d'une courbe $\gamma$ dans la d\'efinition \ref{uple} sont des entiers naturels.

+ D\'efinir ``deux courbes \'equivalentes'' :
 Soit $S_\nu$ une sous-vari\'et\'e de $S_F$ et $\kappa = (i_1, \ldots, i_p)[j_1, \ldots, j_q]$ une fa{\c c}on de $S_\nu$.
Supposons    $\gamma, \gamma': (0, +\infty) \, \to \C^n_{(x)}$ deux courbes tendant vers l'infini avec la m\^eme fa\c con 
$\kappa$  telles que $F \circ \gamma$ et $F \circ \gamma'$ tendent vers deux points de $S_\nu$. 
 Notons  $(l_{i_1}, \ldots, l_{i_p}, l_{j_1}, \ldots, l_{j_q})$ et $(l'_{i_1}, \ldots, l'_{i_p}, l'_{j_1}, \ldots, l'_{j_q})$ leurs  deux $(p+q)$-uples associ\'es aux courbes $\gamma$ et $\gamma'$, respectivement. 
Nous disons que les deux courbes $\gamma$ et $\gamma'$ sont { \it  \'equivalentes}  si leur deux   $(p+q)$-uples associ\'es sont proportionnels, c'est-\`a-dire
\begin{equationth} \label{s:pquple}
\gamma \sim \gamma' \Leftrightarrow (l_{i_1}, \ldots, l_{i_p}, l_{j_1}, \ldots, l_{j_q}) = \lambda (l'_{i_1}, \ldots, l'_{i_p}, l'_{j_1}, \ldots, l'_{j_q}), \quad \text{ o\`u } \lambda \in \C \setminus \{ 0 \} 
\end{equationth}   
(D\'efinition 2.4 de \cite{Thuy}). 

+ D\'efinir ``un point g\'en\'erique d'une sous-vari\'et\'e'' de $S_F$ : Soient $S_\nu$ une sous-vari\'et\'e de $S_F$  et $\kappa$ une fa{\c c}on de $S_\nu$. Soit $a$ un point de $S_\nu$ et $\gamma : (0, +\infty) \, \to \C^n_{(x)}$ une courbe tendant vers l'infini avec la fa{\c c}on $\kappa$ telle que $F \circ \gamma$ tende vers $a$. 
Soit $\{ \hat{\gamma}_i \}$ l'ensemble des courbes \'equivalentes \`a $\gamma$ 
Nous disons que le point $a$ est 
{\it un point g\'en\'erique de  $S_\nu$ relativement \`a la fa{\c c}on $\kappa$} 
s’il admet un voisinage ouvert $U_a$ de $a$ dans $S_{\nu}$ 
qui est un ensemble de limites de courbes $F \circ \hat{\gamma}$, 
o\`u $\hat{\gamma}$ 
est \'equivalente \`a $\gamma$ 
(D\'efinition 2.7 de \cite{Thuy}). 

L'ensemble des points g\'en\'eriques de  $S_\nu$ relativement \`a une fa{\c c}on $\kappa$ 
 est dense dans $S_\nu$.

+ Formaliser l'id\'ee des ``fa{\c c}ons \'etoile'', comme suite: L'id\'ee de la Proposition \ref{remarkfaconetoile}  est la suivante~: \'Etant donn\'ee
  une fa{\c c}on fix\'ee $\kappa$ d'un \'el\'ement $S_\nu$ de la partition de $S_F$ d\'efinie par la relation (\ref{partition facon}), nous subdivisons  $S_\nu$ en utilisant la relation d'\'equivalence  (\ref{s:pquple}) entre  les courbes $\gamma: (0, + \infty) \rightarrow \C^n_{(x)}$ admettant la fa{\c c}on $\kappa$. L'image (par l'application $F$) de 
chaque classe d'\'equivalence   de la relation (\ref{s:pquple})   d\'etermine une sous-vari\'et\'e de $S_\nu$, 
laquelle a une dimension d'autant plus petite que la vitesse de tendre vers z\'ero des courbes coordonn\'ees est grande. Nous obtenons ainsi une 
 sous-partition de  
 $S_{\nu}$. 
\`A partir de ces sous-partitions, nous avons une nouvelle, ``bonne'' partition de $S_F$ (D\'efinition \ref{defXi*}). C'est-\`a-dire, cette fois, cette partition est une stratification (Th\'eor\`eme \ref{theostraXi*}). Le proc\'ed\'e de classer les courbes tendant vers l'infini dans la source, qui sera formalis\'e sous forme de ``fa{\c c}ons \'etoile'' (D\'efinition  \ref{defpreetoile}) est tr\`es significatif~: chaque classe d'\'equivalence de la relation  (\ref{s:pquple}) contient des courbes parall\`eles localement, qui autrement dit, d\'efinissent un feuilletage de dimension (complexe) 1 dans $\C^n_{(x)}$ ({\it cf.}  Exemple \ref{exfeuilletage}). Ce fait est la cl\'e de la d\'emonstration du r\'esultat principal de cet article~: le Th\'eor\`eme \ref{theostraXi*}.

Proc\'ed\'e: Consid\'erons $S_\nu$ un \'el\'ement de dimension $\nu$ de la partition de $S_F$ d\'efinie par la relation (\ref{partition facon}) 
et consid\'erons une fa{\c c}on $\kappa$ de $S_\nu$. L'ensemble des points g\'en\'eriques de $S_\nu$  relativement \`a la fa{\c c}on $\kappa$ 
est dense dans $S_\nu$. Nous avons 
 $$ \dim S^{\kappa}_{\nu_0} = \dim S_{\nu}, \quad \text{ et } \quad \overline{S^{\kappa}_{\nu_0}} = S_{\nu}.$$
Si $S^{\kappa}_{\nu_0} = S_{\nu}$, alors la d\'emonstration est finie. Sinon, nous r\'ep\'etons le proc\'ed\'e ci-dessus~: Notons $S^{\kappa}_{\nu_1}$ l'ensemble des points g\'en\'eriques de $S_{\nu} \setminus S^{\kappa}_{\nu_0}$. Nous avons 
 $$ \dim S^{\kappa}_{\nu_1} = \dim (S_{\nu} \setminus S^{\kappa}_{\nu_0}), \quad \text{ et } \quad \overline{S^{\kappa}_{\nu_1}} = S_{\nu} \setminus S^{\kappa}_{\nu_0}.$$
De plus, nous avons 
$$\quad S^{\kappa}_{\nu_0} \cap S^{\kappa}_{\nu_1} = \emptyset \quad \text{ et } \quad S^{\kappa}_{\nu_1} \subset \overline{S^{\kappa}_{\nu_0}}.$$
Puisque $S^{\kappa}_{\nu_0}$ est dense dans $S_{\nu}$, il vient $\dim (S_{\nu} \setminus S^{\kappa}_{\nu_0}) < \dim S^{\kappa}_{\nu_0}$, nous avons alors 
$$\dim S^{\kappa}_{\nu_0} > \dim S^{\kappa}_{\nu_1}.$$
Si $S^{\kappa}_{\nu_1} = S^{\kappa}_{\nu_0}$, la d\'emonstration est finie. Sinon, nous continuons ce proc\'ed\'e. Puisque $\nu$ est fini,   il existe $t \leq \nu$, $\, t \in \N \setminus \{ 0 \}$ tel que  
les points g\'en\'eriques de  $S^{\kappa}_{\nu_t} $ est  $S^{\kappa}_{\nu_t} $. 
 Donc ${\{ S^{\kappa}_{\nu_i}}\}_{i = 0, \ldots ,t}$ d\'ecrit une partition finie de $S_\nu$. 

 Elle nous permettra de d\'efinir une sous-partition de la partition de $S_F$ d\'efine par les fa{\c c}ons, laquelle sous-partition se r\'ev\'elera \^etre une ``bonne'' stratification.

\begin{proposition} \cite{Thuy} \label{remarkfaconetoile}
{\rm Soit $S_\nu$ un \'el\'ement de dimension $\nu$ de la partition de $S_F$ d\'efinie par la relation (\ref{partition facon}). Pour chaque fa{\c c}on $\kappa$ de $S_\nu$, la relation d'\'equivalence  (\ref{s:pquple}) 
 nous fournit une partition finie ${\{ S^{\kappa}_{\nu_i}}\}_{i = 0, \ldots, t}$ de $S_\nu$,  o\`u $t \leq \nu$, $\, t \in \N$, telle que 
\begin{enumerate}

\item[1)] $\nu = \dim S^{\kappa}_{\nu_0} > \dim S^{\kappa}_{\nu_1} > \dim S^{\kappa}_{\nu_2} > \cdots > \dim S^{\kappa}_{\nu_t},$ 

\item[2)]$S^{\kappa}_{\nu_i} \cap S^{\kappa}_{\nu_j} = \emptyset, \text{ pour }  0 \leq i, j \leq t   \text{ et }  i \ne j, $

\item[3)] $S^{\kappa}_{\nu_i} \subset \overline{S^{\kappa}_{\nu_j}} \text{ pour } i > j \text{ et } 0 \leq i, j \leq t.$
\end{enumerate}

La partition ${\{ S^{\kappa}_{\nu_i}}\}_{i = 0, \ldots, t}$ est appel\'ee
 {\it partition de $S_\nu$ d\'efinie par la fa{\c c}on $\kappa$}.

}
\end{proposition}

\begin{definition} \cite{Thuy} \label{defpreetoile}
{\rm Soit $S_\nu$ un \'el\'ement de la partition de $S_F$ d\'efinie par la relation (\ref{partition facon})  et soit  ${\{ S^{\kappa}_{\nu_i}}\}_{i = 0, ...,t}$  la partition de $S_\nu$ d\'efinie par une fa{\c c}on $\kappa$ de $S_\nu$ comme dans la Proposition \ref{remarkfaconetoile}. Si $t \geq 1$, 
nous d\'efinissons  les ``{\it fa{\c c}ons \'etoile}'' de la fa{\c c}on $\kappa$, comme suit~: 

la fa{\c c}on $\kappa$ de $ S^{\kappa}_{\nu_1}$ appel\'ee la {\it fa{\c c}on \'etoile} $\kappa^{1*}$ de $S_\nu$,

la fa{\c c}on de $\kappa$ de $ S^{\kappa}_{\nu_2}$ appel\'ee la {\it fa{\c c}on \'etoile} $\kappa^{2*}$ de $S_\nu$,

$\quad \quad \quad \quad \quad \quad \quad \quad$ ...

la fa{\c c}on $\kappa$ de $ S^{\kappa}_{\nu_t}$ appel\'ee la {\it fa{\c c}on \'etoile} $\kappa^{t*}$ de $S_\nu$. 

Par convention, nous disons que la fa{\c c}on $\kappa$ est la fa{\c c}ons $\kappa^{0*}$ de $S_\nu$.
}
\end{definition}

\begin{example} \label{exemplesubdivisionetoile1}
{\rm La fa{\c c}on (2)[1] de l'ensemble asymptotique de l'application $F(x_1,x_2) = \left({(x_1x_2)}^2, {(x_1x_2)}^3 + x_1 \right)$  
  admet une seule fa{\c c}on \'etoile ${(2)[1]}^{1*}$.
}
\end{example}

Chaque \'el\'ement $S_{\nu_i}^{\kappa}$ de la partition ${S_{\nu_i}^{\kappa}}_{i = 0, \ldots, t}$
  d\'efini dans la Proposition \ref{remarkfaconetoile}
 peut \^etre associ\'e \`a plusieurs classes d’\'equivalence de courbes, 
relativement \`a la relation d'\'equivalence \label{s:pquple}, 
mais nous fixons seulement une classe d'\'equivalence de courbes correspondantes. 
Avec cette convention, nous avons :

\begin{proposition} \cite{Thuy} \label{feuilletage}
{\rm Fixons une fa{\c c}on \'etoile $\kappa^{i*}$ de $S_F$. 
Soit $U_a$ un voisinage ouvert suf\-fisamment petit d'un point $a \in S_F$ 
tel que tout point de $U_a$ admet la fa{\c c}on \'etoile $\kappa^{i*}$. 
Alors toutes les courbes correspondantes aux points de $U_a$ d\'efinissent localement un feuilletage de dimension (complexe) 1 de l'espace source $\C^n_{(x)}$.

}
\end{proposition}

\begin{definition} \label{defXi*} 
{\rm Consid\'erons  $S_\nu$ un \'el\'ement de la stratification de $S_F$ d\'efinie par les fa{\c c}ons
et $\Xi(S_\nu)$ l'ensemble  de toutes les fa{\c c}ons de $S_\nu$. 
Nous d\'efinissons $S_{\nu_0}$ comme l'ensemble des points g\'en\'eriques relativement \`a  
au moins une fa{\c c}on de $S_\nu$, c'est-\`a-dire ~: 
$$S_{\nu_0} = \bigcup_{\kappa \, \in \, \Xi(S_\nu)} S^{\kappa}_{\nu_0},$$
 o\`u $S_{\nu_0}^\kappa$ est d\'efini comme dans la Proposition \ref{remarkfaconetoile}.

 Notons $A_{\nu_0} = S_\nu \setminus S_{\nu_0}$. 
Nous d\'efinissons $S_{\nu_1}$ comme l'ensemble des points g\'en\'eriques  relativement \`a au moins une fa{\c c}on de $A_{\nu_0}$ 
et nous notons $A_{\nu_1} = A_{\nu_0} \setminus S_{\nu_1}$. En g\'en\'eral, pour $i \geq 1$, 
nous d\'efinissons 
$S_{\nu_i}$ comme l'ensemble des points g\'en\'eriques  relativement \`a au moins une fa{\c c}on de $A_{\nu_{i-1}} $ et 
$A_{\nu_i} = A_{\nu_{i-1}} \setminus S_{\nu_i}$. Nous obtenons une   sous-partition de $S_\nu$. 
 \`A partir de ces sous-partitions, nous obtenons ainsi une nouvelle partition de $S_F$, 
appel\'ee {\it la partition de $S_F$ d\'efinie par les fa{\c c}ons \'etoile}. 
}
\end{definition}

\begin{theorem} \label{theostraXi*}  \cite{Thuy}
Soit $F : \C_{(x)}^n \to \C_{(\alpha)}^n$  une application polynomiale dominante. La partition de l'ensemble asymptotique de $F$ d\'efinie par les fa{\c c}ons \'etoile est une stratification  diff\'erentiable satisfaisant la propri\'et\'e de fronti\`ere.
\end{theorem}

\section{Une stratification de Thom-Mather de l'ensemble asymptotique}

Consid\'erons  la stratification $(\mathscr{S})$  de $S_F$ d\'efinie par les fa{\c c}ons \'etoile. Nous allons prouver que la stratification $(\mathscr{S})$ est une stratification  de Thom-Mather. 
Pour faire cela, nous devons construire des voisinages tubulaires des strates et les fonctions $\rho$ et $\pi$ 
correspondantes (cf. D\'efinition \ref{definitionthommather}).
 En fait, nous allons construire les rayons des voisinages tubulaires et les fonctions 
``distance'' $\rho$ corrrespondant aux startes. 
L'id\'ee est d'utiliser les courbes correspondantes aux fa{\c c}ons 
(Lemme \ref{lemmecourbe}). Les courbes correspondantes de chaque strate de $(\mathscr{S})$ forment un feuilletage (Proposition  \ref{feuilletage}). 
L'image  de ce feuilletage fournira les rayons des vosinages tubulaires. Les fonctions ``distance''  $\rho$ sont
d\'efinies par r\'ecurrence. Les fonction ``projection'' $\pi$ sont d\'efinies automatiquement le long des rayons des vosinages tubulaires, ou, autrement dit le long des images des coubres correspondantes aux fa{\c c}ons. 
Les voisinages obtenus sont du type des voisinages tubulaires effil\'es d\'efinis 
par Marie-H\'el\`ene Schwartz dans \cite{MHS, MHS2}. Pour formaliser cette id\'ee, nous avons besoin d'une relation entre 
des strates adjacentes. Nous d\'efinissons une relation d'ordre partielle entre les fa{\c c}ons des strates de $(\mathscr{S})$ dans la partie suivante.

\subsection{La relation entre les fa{\c c}ons des strates de la stratification de $S_F$ d\'efinie par les fa{\c c}ons} \label{relationd'ordre}

Notons que dans la suite, quand nous concernons \`a une fa{\c c}on $\kappa$, alors $\kappa$ peut \^etre une fa{\c c}on \'etoile, au sens de la D\'efinition \ref{defpreetoile}.

\begin{definition} \label{defprec}
{\rm Soient $\kappa, \kappa' $ deux fa{\c c}ons de $S_F$ telles que $\kappa = (I_p)[J_q],$ et $\kappa' = (I_{p'})[J_{q'}].$ 
D\'efinissons la relation d'ordre partielle $\kappa \prec \kappa'$ si nous  avons l'un des trois cas suivants :

1) $\{ I_p\} \supsetneq \{ I_{p'} \}$ et $\{ J_q\} \supset \{ J_{q'} \}.$

2) $\{ I_p\} = \{ I_{p'} \}$ et $\{ J_q\} \supsetneq \{ J_{q'} \}.$ 

3) $\kappa$ est une fa{\c c}on \'etoile de $\kappa'$.
}
\end{definition}

\begin{definition} \label{defenitionprec}
{\rm Nous disons que $a \prec a'$, o\`u $a, a' \in S_F$, si pour tout $\kappa' \in \Xi^*(a')$, il existe $\kappa \in \Xi^*(a)$ tel que $ \kappa = \kappa'$ ou $ \kappa \prec \kappa'$.}
\end{definition}

\begin{definition}
{\rm Nous d\'efinissons {\it l'ordre}  du point $a$ de $S_F$, notons $or(a)$, 
comme le nombre des fa{\c c}ons de $a$.
}
\end{definition}

\begin{remark} \label{reprec}
{\rm D'apr\`es la D\'efinition \ref{defprec}, alors si $a \prec a'$, nous avons $or(a) \geq or(a')$.}
\end{remark}

\begin{theorem} \label{propofrontiere}
Soient $a, a' \in S_F$. Alors $a \prec a'$ si et seulement si $strate(a) \subset \overline{strate(a')},$
o\`u $strate(a)$ et $strate(a')$ sont des strates de la stratification de $S_F$ d\'efinie par les fa{\c c}ons \'etoile, contiennent $a$ et $a'$, respectivement. 
\end{theorem}

\begin{preuve}
Soient $a, a' \in S_F$ tels que $a \prec a'$. Nous avons les deux cas suivants~:

1) $or(a) > or (a')$ : La Propositon \ref{proalgebrique} dit que chaque fa{\c c}on de $S_F$ d\'etermine une \'equation alg\'ebrique obtenue \`a partir de facteurs de l'\'equation de l'ensemble $S_F$ dans le Th\'eor\`eme \ref{theoremjelonek1}. Comme $or(a) > or (a')$, et par la construction des strates de la stratification de $S_F$ d\'efinie par les fa{\c c}ons \'etoile (D\'efinition \ref{defXi*}), l'ensemble des \'equations du syst\`eme d'\'equations de $\overline{strate(a)}$ contient  l'ensemble des \'equations du syst\`eme des \'equations de $\overline{strate(a')}$. 
Cela signifie que
$strate(a) \subset \overline{strate(a')}.$

2) $or(a) = or(a')$ : il existe une fa{\c c}on $\kappa = (I_p)(J_q)$ de $strate(a)$ et une fa{\c c}on $\kappa' = (I'_p)(J'_q)$ de $strate(a')$ telles que $\kappa \prec \kappa'.$ Dans ce cas-ci, par la D\'efinition \ref{defprec}, nous avons deux possibilit\'es~:

a) ou bien $I_p \cup J_q \supsetneq I'_p \cup J'_q$ : la fa{\c c}on $\kappa'$ n'est pas une fa{\c c}on \'etoile de $\kappa$. 
 Il existe $a'' \in S_F$ tel qu'une fa{\c c}on 
$\kappa'' = (I^{''}_p) [J^{''}_q]$ de $a''$ 
satisfait $I'_p \cup I''_p = I_p$, $J'_q \cup J''_q = J_q$ (unions non n\'ecessairement disjointes). 
 Cela signifie que $\overline{strate(a)} = \overline{strate(a')} \cap \overline{strate(a'')}$.  Nous avons donc  $strate(a) \subset \overline{strate(a')}.$

b) ou bien $I_p \cup J_q = I'_p \cup J'_q$ : la fa{\c c}on $\kappa$ est une fa{\c c}on  \'etoile de $\kappa'$.  Par la D\'efinition \ref{defpreetoile} e la D\'efinition \ref{defXi*}, la strate $strate(a)$ est incluse dans $\overline{strate(a')}.$
\end{preuve}

Avant d'utiliser la relation entre les fa{\c c}ons dans le Th\'eor\`eme \ref{propofrontiere} ci-dessus pour contruire des voisinages tubulaires des strates de la stratification de $S_F$ d\'efinie par les fa{\c c}ons \'etoile, nous avons besoin de deux lemmes suivants. 

\subsection{Lemmes}
\begin{lemma} \label{lemmethuy2}
Soit $F: \C^n_{(x)} \to \C^n_{(\alpha)}$ une application polynomiale dominante, 
alors $$  S_F \cup F(\C^n_{(x)}) = \C^n_{(\alpha)}.$$
\end{lemma}
 
\begin{preuve}
Puisque $F$ est dominante, nous avons  $\overline{F(\C^n_{(x)})} = \C^n_{(\alpha)}$ (cf. D\'efinition \ref{definitiondominant}). 
Prenons $a \in \overline{F(\C^n_{(x)})} \setminus F(\C^n_{(x)})$, 
alors il existe une suite $ \{\xi_k\} \subset \C^n_{(x)}$ 
telle que $F(\xi_k)$ tend vers $a$. 
Si $ \xi_k $ ne tend pas vers l'infini, 
alors $ \xi_k $ tend vers $x_0 \in \C^n_{(x)}$
et  $F(\xi_k)$ tend vers $F(x_0)$. 
Nous avons donc $F(x_0) = a$, d'o\`u la contradiction avec $a \notin F(\C^n_{(x)})$. 
Donc $ \xi_k $ tend vers l'infini et $a \in S_F$. 
\end{preuve}

\begin{lemma} \label{lemmethuy3}
Soit $a\in \C^n_{(\alpha)} \setminus S_F$, soit $\kappa_\nu$ une fa\c con d'une strate $S_\nu$ de la stratification de $S_F$ d\'efinie par les fa{\c c}ons \'etoile et soit
$a^\nu$ un point de $S_\nu$, suffisamment proche de $a$. Alors, nous pouvons toujours d\'eterminer une courbe diff\'erentiable 
$$\gamma^{\nu} : [1, + \infty ) \to \C^n_{(x)}, $$
telle que, d'une part, lorsque $u$ tend vers l'infini, $\gamma^{\nu}(u)$ tend vers l'infini avec la fa{\c c}on $\kappa_{\nu}$ et $F(\gamma^{\nu}(u))$ tend vers $a^\nu$, et d'autre part,  $F({\gamma}^{\nu}(1)) = a$.
\end{lemma}

\begin{preuve}
D'abord, puisque $a^\nu \in S_F$, alors, par le Lemme \ref{lemmecourbe}, Il existe une courbe diff\'erentiable 
$$\gamma^{\nu} : [1, + \infty ) \to \C^n_{(x)}, $$
telle que  $\gamma^{\nu}(u)$ tend vers l'infini avec la fa{\c c}on $\kappa_{\nu}$ et $F(\gamma^{\nu}(u))$ tend vers $a^\nu$ lorsque $u$ tend vers l'infini. Comme $a\in \C^n_{(\alpha)} \setminus S_F$, alors par le Lemme \ref{lemmethuy2}, nous avons $a \in F(\C^n_{(x)})$. D\'eterminons  maintenant la courbe $\gamma^{\nu}$ telle que ${\gamma}^{\nu}(1)$ est contenu dans $F^{-1}(a)$, c'est-\`a-dire,  $F({\gamma}^{\nu}(1)) = a$. 
 D'apr\`es le Th\'eor\`eme \ref{theostraXi*},  $\overline{S_\nu}$ est une vari\'et\'e alg\'ebrique. Nous pouvons supposer que l'\'equation de $\overline{S_\nu}$ est 
$$\phi(\alpha) = 0$$ 
(notons que, par ``l'\'equation de $\overline{S_\nu}$'', nous entendons possiblement un syst\`eme d'\'equations). Sans perte de g\'en\'eralit\'e, nous pouvons supposer que  
$$\kappa_\nu = (1, \ldots, p)[p+1, \ldots, p+q],$$
o\`u $q+p \leq n.$ La courbe $\gamma^{\nu}$ peut s'\'ecrire sous la forme :
{
\begin{equationth} \label{formecourbe7}
\begin{split}
( f_{p+1}(a^{\nu}) + \frac{1}{{\lambda^{{\nu}}_1} u}, \ldots, f_{p+q}(a^{\nu}) + \frac{1}{{\lambda^{{\nu}}_q}u},
\hskip 5.5truecm  \\
f_{1}(a^{\nu}, \lambda^{{\nu}}_1, \ldots, \lambda^{{\nu}}_q, u), \ldots,  f_{p}(a^{\nu}, \lambda^{{\nu}}_1, \ldots,  \lambda^{{\nu}}_q, u), f_{p+q+1}, \ldots, f_n ),
\end{split}
\end{equationth}
}
\noindent o\`u 
$f_{l}(a^\nu, \lambda^{\nu}_1, \ldots,  \lambda^{\nu}_q, u)$ tend vers l'infini  pour $l = 1, \ldots, p$ et $f_{p+q+1}, \ldots, f_n$ tendent vers des nombres complexes fix\'es lorsque $u$ tend vers l'infini. 
Puisque $F({\gamma}^{\nu}(1)) = a$ et 
$a^\nu \in S_{\nu}$ alors $(\lambda^{\nu}_1, \ldots,  \lambda^{\nu}_q)$ est la solution du syst\`eme d'\'equations 
\begin{equation*}
\begin{cases}
F({\gamma}^{\nu}(1)) = a \cr
\phi(a^\nu) = 0
\end{cases},
\end{equation*}
nous pouvons donc toujours d\'eterminer la courbe ${\gamma}^{\nu}$.

\end{preuve}

L'exemple suivant illustre le Lemme \ref{lemmethuy3} ci-dessus. 

\begin{example} \label{exeThomMather1}
{\rm Consid\'erons l'application polynomiale dominante $F : \C^3_{(x_1,x_2,x_3)} \to \C^3_{(\alpha_1, \alpha_2, \alpha_3)}$ telle que $$F(x_1,x_2,x_3) = (x_1^2 - 1, x_2+2, (x_1^2 -1)(x_2+2)x_3).$$ 
L'ensemble asymptotique $S_F$ de $F$ est l'union de deux plans $\{ \alpha_1= 0\}$ et $\{ \alpha_2 = 0\}$. Prenons le point $a^1 = (0, 0, 2)$ appartient \`a une strate de dimension 1 de $S_F$ et nous voyons que $a^1 \in  \{ \alpha_1= 0\} \cap \{ \alpha_2 = 0\}$. Nous pouvons v\'erifier facilement qu'une fa{\c c}on de $a^1$ est $\kappa = (3)[1, 2]$. Prenons le point $a = (3, 1, 3) \in \C^3_{(\alpha_1, \alpha_2, \alpha_3)} \setminus S_F$. 
D'apr\`es le lemme \ref{lemmethuy2}, le point $a$ est contenu dans $F(\C^3_{(x_1,x_2,x_3)})$. 
Cherchons maintenant une courbe ${\gamma}^1 : [1, + \infty ) \to \C^3_{(x_1, x_2, x_3)}$  
 telle que, d'une part, ${\gamma}^1(u)$ tende vers l'infini avec la fa{\c c}on $(3)[1,2]$ et $F({\gamma}^1(u))$ tende vers le point $a^1$  lorsque $u$ tend vers l'infini, et que, d'autre part, $F({\gamma}^1(1))$ soit $a$. Comme ${\gamma}^1(u)$ tende vers l'infini avec la fa{\c c}on $(3)[1,2]$ et  les deux premi\`eres coordonn\'ees de $F({\gamma}^1(u))$ tendent vers 0, donc ${\gamma}^1(u)$ peut \^etre \'ecrit sous la forme: 
$$ {\gamma}^1(u) = \left( 1 + \frac{1}{\lambda u}, -2 + \frac{1}{\mu u}, x_{3,u} \right), \, \text{ o\`u } x_{3,u} \to \infty \quad \text{ et } \lambda, \mu \in \C.$$
Nous avons 
$$ F({\gamma}^1(u)) = \left( \frac{2}{\lambda u} + \frac{1}{\lambda^2 u^2}, \frac{1}{\mu u}, \left( \frac{2}{\lambda u} + \frac{1}{\lambda^2 u^2} \right) \frac{1}{\mu u} x_{3,u} \right).$$ 
Puisque $F({\gamma}^1(u))$ tend vers $a^1 = (0, 0, 2)$, il vient~:
$$ \begin{aligned}
& x_{3,u} = \lambda \mu u^2, \cr
& F({\gamma}^1(u)) = \left( \frac{2}{\lambda u} + \frac{1}{\lambda^2 u^2}, \frac{1}{\mu u}, 2 + \frac{1}{ \lambda u} \right). 
\end{aligned}$$
Puisque $F({\gamma}^1(1)) = a = (3, 1, 3)$,  nous avons $\lambda = \mu = 1$. Nous avons donc  $$\gamma^1(u) = \left(1 + \frac{1}{u}, -2 +\frac{1}{u}, u^2 \right).$$
}
\end{example}

\subsection{Th\'eor\`eme sur une stratification de Thom-Mather de l'ensemble asymptotique} \label{preuveThomMather} 
\begin{theorem} \label{theoremThom-Mather}
Soit $F: \C^n_{(x)} \to \C^n_{(\alpha)}$ une application polynomiale dominante. La stra\-ti\-fi\-ca\-tion de de l'ensemble asymptotique de $F$  
d\'efinie par les fa{\c c}ons \'etoile est une stratification de Thom-Mather. 
\end{theorem}

\begin{preuve}

La premi\`ere \'etape consiste d'une part  \`a attacher \`a chaque strate $S_\nu$ une ou plusieurs fa\c cons de $S_\nu$
de mani\`ere \`a ce que pour chaque s\'equence de strates incluant $S_\nu$
$$\emptyset \subset S_0 \subset \overline{S_1} \subset \cdots \subset \overline{S_\nu} \subset \cdots \subset \overline{S_{n-1}}$$
corresponde une suite de fa\c cons 
$$\kappa_0 \prec \kappa_1 \prec \cdots \prec \kappa_{n-1}.$$
D'autre part, nous consid\'erons des voisinages $V_\nu$ 
des strates $S_\nu$ de $S_F$, effil\'es au sens de M.-H. Schwartz (\cite{MHS2}). Les voisinages 
tubulaires que nous allons construire en sont des ``sous-tubes". 

La d\'emonstration est faite par r\'ecurrence
d\'ecroissante, en commen\c cant par la strate de dimension la plus grande : $S_{n-1}$, que nous notons $S_{n_0}$. 
A chaque \'etape, nous d\'efinissons le voisinage tubulaire ainsi que les fonctions $\rho$ et $\pi$. 

Consid\'erons donc un point $a$ situ\'e dans $V_{n_0}\setminus (\bigcup_{\nu < n-1} V_\nu \cup S_F)$. 
Nous pouvons choisir un point $a^{n_0}$ suffisamment proche de $a$ et, 
pour toute fa\c con $\kappa_{n_0}$ de $S_{n-1}= S_{n_0}$, le Lemme \ref{lemmethuy3} nous 
fournit le rayon $F(\gamma^{n_0}(u))$  reliant $a$ \`a $a^{n_0}$. Le voisinage tubulaire ${\mathcal T}_{n_0}$ de $S_{n_0}\setminus \bigcup_{\nu < n-1} V_\nu$ est constitu\'e de 
l'ensemble des rayons $F(\gamma^{n_0}(u))$
aboutissant aux points $a^{n_0}$ de $S_{n_0}\setminus \bigcup_{\nu < n-1} V_\nu$. 
Pour tout point $a' = F(\gamma^{n_0}(u))$ situ\'e sur ce rayon, 
nous d\'efinissons 
\begin{equationth} \label{pirho0} 
\pi_{n_0}(a') = a^{n_0} \text{ et } \rho_{n_0}(a') = d_c(a', a^{n_0})
\end{equationth} 
o\`u $d_c(a', a^{n_0})$ est la distance curviligne, de $a'$ \`a $a^{n_0}$, le long de la courbe $F(\gamma^{n_0})$. 

Le cas g\'en\'eral de la r\'ecurrence suit les m\^emes techniques que dans le cas de deux strates. Celui-ci est 
cependant plus simple \`a expliciter, ce que nous faisons ci-dessous. Nous examinons ensuite le cas g\'en\'eral de la 
r\'ecurrence.

Consid\'erons donc maintenant le cas de deux strates, plus pr\'ecis\'ement, la strate $S_{n_0}$ et une strate $S_{n_1}$ de 
dimension imm\'ediatement inf\'erieure \`a celle de $S_{n_0}$. 
Nous allons construire le voisinage tubulaire $ {\mathcal T}_{n_1}$ 
de $S_{n_1}$ en trois parties, la premi\`ere \'etant la plus d\'elicate. 

1) Nous notons $\kappa_{{n_0}}$ et $\kappa_{n_1}$ deux fa\c cons de $S_{{n_0}}$ et $S_{n_1}$ respectivement, telles que 
$\kappa_{n_1} \prec \kappa_{{n_0}}.$ 
Consid\'erons un point $a$ situ\'e dans $(V_{{n_0}} \cap V_{n_1})\setminus S_F$. 
Nous pouvons, comme pr\'e\-c\'e\-dem\-ment, choisir un point $a^{{n_0}}$ suffisamment proche de $a$,
et utiliser  le Lemme \ref{lemmethuy3}. Nous pouvons supposer que  
$$\kappa_{{n_0}} = (1, \ldots, p)[p+1, \ldots, p+q],$$
o\`u $q+p \leq n.$ La courbe $\gamma^{{{0}}}$ est construite comme dans la preuve du Lemme \ref{lemmethuy3} et 
peut s'\'ecrire sous la forme :
{
\begin{equationth} \label{formecourbe7}
\begin{split}
( f_{p+1}(a^{n_0}) + \frac{1}{{(\lambda^{{n_0}}_1)}u}, \ldots, f_{p+q}(a^{n_0}) + \frac{1}{{(\lambda^{{n_0}}_q)}u},
\hskip 5.5truecm  \\
f_{1}(a^{n_0}, \lambda^{{n_0}}_1, \ldots, \lambda^{{n_0}}_q, u), \ldots,  f_{p}(a^{n_0}, \lambda^{{n_0}}_1, \ldots,  \lambda^{{n_0}}_q, u), f_{p+q+1}, \ldots, f_n ),
\end{split}
\end{equationth}
}
\noindent o\`u 
$f_{l}(a^{n_0}, \lambda^{{n_0}}_1, \ldots,  \lambda^{{n_0}}_q, u)$ tend vers l'infini  pour $l = 1, \ldots, p$ et $f_{p+q+1}, \ldots, f_n$ tendent vers des nombres complexes fix\'es lorsque $u$ tend vers l'infini. 
Nous obtenons ainsi le rayon $(a, a^{{n_0}})$ du voisinage tubulaire ${\mathcal T}_{{n_0}}$.

Le voisinage tubulaire ${\mathcal T}_{{n_0}}$ de $S_{n_0}$ dans $V_{{n_0}} \cap V_{n_1}$ 
est construit de fa\c con \`a ce que 
la longueur $\varepsilon (a^{n_0})$ 
des rayons de ${\mathcal T}_{{n_0}}$ tende vers $0$ lorsque la distance du  point $a^{n_0}$ \`a la strate $S_{n_1}$  tend vers $0$. Cette construction correspond \`a  celle des voisinages tubulaires effil\'es de
M.-H. Schwartz (\cite{MHS2}).

Appelons ${\mathcal T}_{{n_0}}^\varepsilon$ le voisinage tubulaire ainsi obtenu. Pour tout $ a^{n_0}$,  notons
$\varepsilon'(a^{n_0}) = \varepsilon(a^{n_0}) /2$. Nous obtenons ainsi un voisinage ${\mathcal T}_{{n_0}}^{\varepsilon'}$  
construit par homoth\'etie de rapport $1/2$ (voir figure \ref{thommather1}).

Nous construisons maintenant ${\mathcal T}_{n_1}$  
et d\'efinissons  les fonctions $\pi_{n_0}, \rho_{n_0}$ et 
$\pi_{n_1}, \rho_{n_1}$ dans ${\mathcal T}_{{n_0}}^{\epsilon'}\cap V_{n_1}$. Pour cela, 
consid\'erons  un point $a$ de ${\mathcal T}_{{n_0}}^{\epsilon'}\cap V_{n_1}$. 
Nous pouvons choisir un point $a^{n_1}$ situ\'e sur la strate $S_{n_1}$, suffisamment proche de $a$,
et utiliser  le Lemme \ref{lemmethuy3}. 
Sans perte de g\'en\'eralit\'e, nous pouvons supposer que 
$$\kappa_{n_1} = (1, \ldots, p)[p+1, \ldots, p+q+1].$$
Nous obtenons une courbe $\gamma^{n_1}(t)$, laquelle  peut s'\'ecrire sous la forme :
{
\begin{equationth} \label{formecourbe3} 
\begin{split}
( g_{p+1}(a^{n_1}) + \frac{1}{{(\lambda^{n_1}_1)}t}, \ldots, g_{p+q+1}(a^{n_1}) + \frac{1}{{(\lambda^{n_1}_{q+1})}t}, 
\hskip 5.5truecm  \\
g_{1}(a^{n_1}, \lambda^{n_1}_1, \ldots, \lambda^{n_1}_{q+1}, t), \ldots,  g_{p}(a^{n_1}, \lambda^{n_1}_1, \ldots,  \lambda^{n_1}_{q+1}, t), g_{p+q+2}, \ldots, g_n ),
\end{split}
\end{equationth}
}
\noindent o\`u 
$g_{l}(a^{n_1}, \lambda^{n_1}_1, \ldots,  \lambda^{n_1}_{q+1}, t)$ tend vers l'infini  pour $l = 1, \ldots, p$ et 
$g_{p+q+2}, \ldots, g_n$ tendent vers des nombres complexes fix\'es lorsque $t$ tend vers l'infini. 
Puisque $F({\gamma}^{{n_1}}(1)) = a$ et 
$a^{n_1} \in S_{{n_1}}$ alors $(\lambda^{n_1}_1, \ldots,  \lambda^{n_1}_{q+1})$ est une solution 
du syst\`eme d'\'equations 
\begin{equation*}, 
\begin{cases}
F({\gamma}^{{n_1}}(1)) = a \cr
\phi_1(a^{n_1}) = 0
\end{cases},
\end{equation*}
o\`u $\phi_1$ est l'\'equation de la vari\'et\'e alg\'ebrique $S_{n_1}$. La courbe $F(\gamma^{n_1}(t))$, allant 
de $a$ \`a $a^{n_1}$ est un rayon du voisinage tubulaire ${\mathcal T}_{n_1}$  de $S_{n_1}$ (voir figure \ref{thommather2}). 

Soit $u_j \in (1, +\infty)$ 
une valeur fix\'ee du param\`etre $u$ de la courbe $\gamma^{{{n_0}}}(u)$, notons 
$a_j$ le point $F(\gamma^{{{n_0}}}(u_j))$ situ\'e sur la courbe $(a, a^{{n_0}})$.
De la m\^eme mani\`ere que la courbe $\gamma^{{n_1}}$ (formule (\ref{formecourbe3})) nous 
pouvons construire des courbes 
$F(\gamma^{{n_1}}_j(t))$ reliant $a_j$  \`a $a^{n_1}$   et avec le 
m\^eme param\`etre $t$  (voir figure \ref{thommather2}). 

\begin{figure}[h!]
\centering
\includegraphics[scale=0.7]{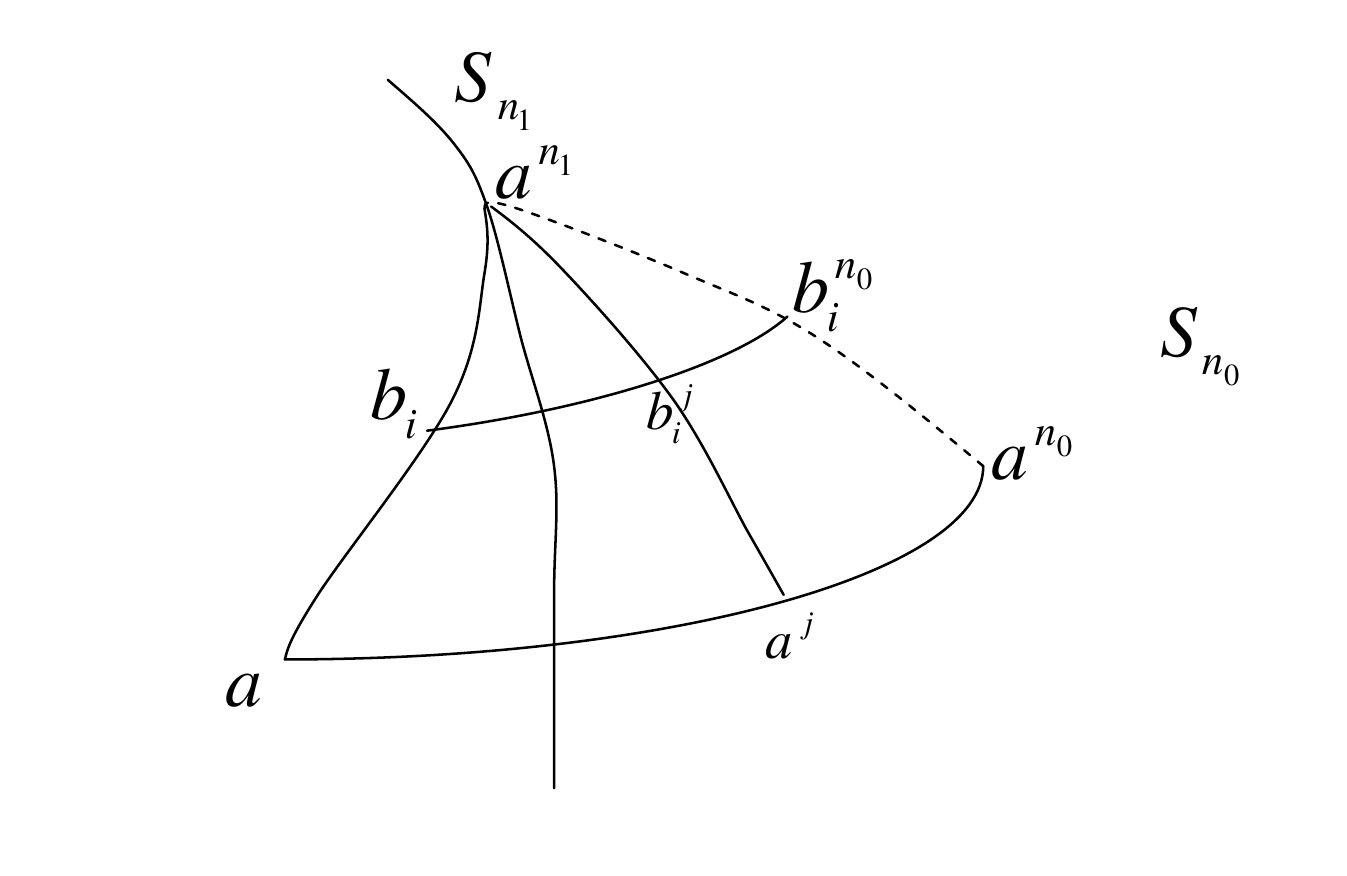}
\caption{Construction de ${\mathcal{T}}_{n_1} \cap {\mathcal{T}}_{n_0}^{\epsilon'}$.}
\label{thommather2}
\end{figure}

D'une part, lorsque $u_j$ 
tend vers l'infini, la courbe $(a_j, a^{n_1})$ tend vers une courbe $(a^{{n_0}}, a^{n_1})$ situ\'ee dans la strate $S_{{n_0}}$. 
Les courbes $(a_j, a^{n_1})$ et leur limite $(a^{{n_0}}, a^{n_1})$ sont des rayons du voisinage tubulaire 
${\mathcal T}_{n_1}$ de la strate $S_{n_1}$, que nous construisons. 
Par d\'efinition, pour tout point $a'$ situ\'e sur une courbe
 $(a_j, a^{n_1})$ (ou sur la limite $(a^{{n_0}}, a^{n_1})$), nous posons $\pi_{n_1} (a') = a^{n_1}$.
 

 D'autre part, tout point $b_i$ situ\'e sur la courbe
 $(a, a^{n_1})$ s'\'ecrit $b_i = F({\gamma}^{n_1}(t_i))$ pour une valeur $t_i$. Fixons cette valeur $t_i$ du 
 param\`etre $t$. Pour $u_j$ allant de $1$ \`a $+\infty$, 
 l'ensemble des points $b_i^j= F({\gamma}^{{n_1}}_j(t_i))$ d\'ecrit une courbe 
  reliant le point $b_i$ \`a un point $b_i^{{n_0}}$ de $S_{n_0}$. Le choix des courbes $ {\gamma}^{{n_0}}$ et
  ${\gamma}^{n_1}$ ((\ref{formecourbe7}) et (\ref{formecourbe3})) nous permet d'affirmer que le point $b_i^{{n_0}} = \lim_{u_j \to \infty} (F({\gamma}^{{n_1}}_j(t_i)))$ 
  est situ\'e sur la courbe $(a^{{n_0}}, a^{n_1})$. Par d\'efinition, nous disons que les points de la courbe 
  $(b_i, b_i^{{n_0}})$ 
 sont situ\'es \`a ``m\^eme distance" du point $a^{n_1}$, relativement \`a  la fonction $\rho_{n_1}$, plus pr\'ecis\'ement nous posons :
 $$\rho_{n_1} (b_i) = \rho_{n_1}(b_i^j) = \rho_{n_1} (b_i^{{n_0}}) = \wideparen{b_i^{n_0} a^{n_1}} \qquad \text{ pour tout } t_i,$$
o\`u $\wideparen{b_i^{n_0} a^{n_1}}$ est le longueur de la courbe reliant deux points $b_i^{n_0}$ et $a^{n_1}$, d\'efini comme le longueur de deux points dans l'espace $\C^n$, avec le sens classique. 
 En particulier, 
$$\rho_{n_1} (a) = \rho_{n_1}(a_j) = \rho_{n_1} (a^{{n_0}}) = \wideparen{a^{n_0} a^{n_1}}.$$

2) Construction de $ {\mathcal T}_{n_1}\setminus  {\mathcal T}_{{n_0}}^{\epsilon}$. 

Dans $V_{n_1} \setminus ({\mathcal T}_{{n_0}}^{\epsilon}\cap V_{n_1} )$, c'est-\`a-dire dans la partie qui ne rencontre
pas les  voisinages tubulaires des strates dont $S_{n_1}$ est adjacente, 
la situation est celle du 
Lemme  \ref{lemmethuy3}, relativement \`a la fa\c con $\kappa_{n_1}$. Nous obtenons des courbes
$\gamma^{n_1}$ dont les images par $F$ sont des rayons de ${\mathcal T}_{n_1}$, ce qui d\'efinit la fonction $\pi_{n_1}$. 
La fonction   $\rho_{n_1}$ 
est d\'efinie comme en (\ref{pirho0}). 

 \begin{figure}[h!]
\centering
\includegraphics[scale=0.7]{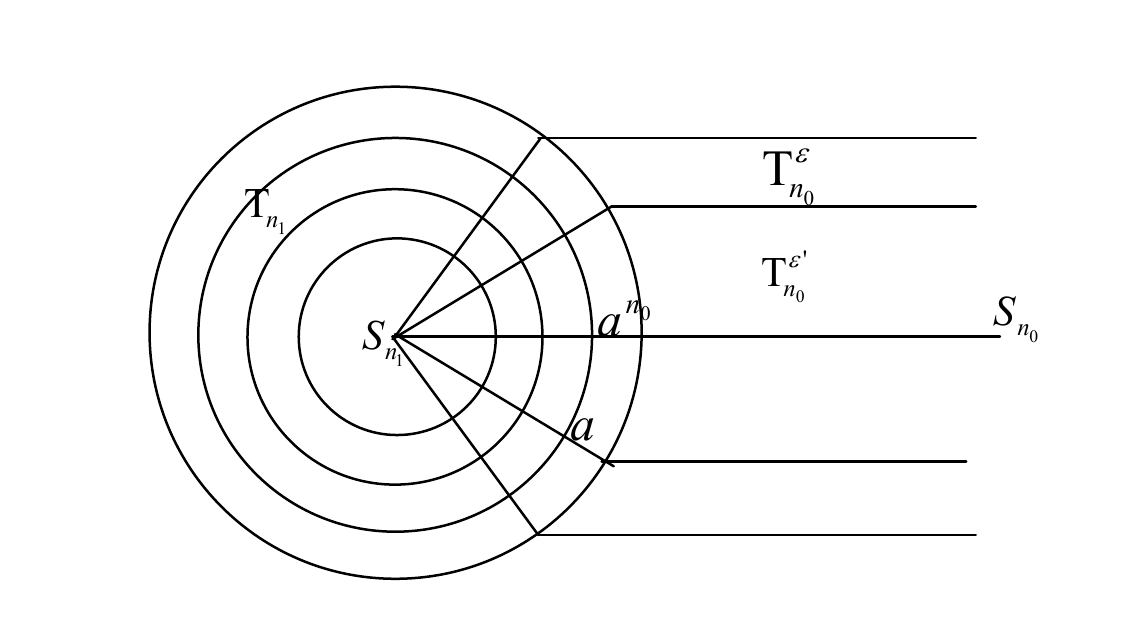}
\caption{Le voisinage tubulaire ${\mathcal T}_{n_1}$ en trois parties.}
\label{thommather1}
\end{figure}

3) Construction de $ {\mathcal T}_{n_1}  \cap ({\mathcal T}_{{n_0}}^{\epsilon} \setminus {\mathcal T}_{{n_0}}^{\epsilon'})$.

De m\^eme que dans le cas pr\'ec\'edent, le Lemme \ref{lemmethuy3} fournit des courbes 
$\gamma^{n_1}$ rayons de ${\mathcal T}_{n_1}$, ce qui d\'efinit la fonction $\pi_{n_1}$.
Dans $V_{n_1} \cap ({\mathcal T}_{{n_0}}^{\epsilon} \setminus {\mathcal T}_{{n_0}}^{\epsilon'})$, la fonction $\rho_{n_1}$ 
est d\'ej\`a d\'efinie sur $V_{n_1} \cap \partial ({\mathcal T}_{{n_0}}^{\epsilon})$ et sur 
$V_{n_1} \cap \partial ({\mathcal T}_{{n_0}}^{\epsilon'})$. Nous la prolongeons
de fa\c con diff\'erentiable, le long des courbes $F(\gamma^{n_0})$ 
en utilisant les fonctions de prolongement de Whitney (\cite{Whitney} IV, \textsection 27)
 (voir figure \ref{thommather1}).
 
 Remarquons que si l'on a deux strates (ou plus) $S_{n_0}$ et $S'_{n_0}$ telles que 
 telles que $S_{n_1} \subset {\overline{S_{n_{0}}}}$ et $S_{n_1} \subset {\overline{S'_{n_0}}}$, alors, 
comme les voisinages $V_{n_0}$ et $V'_{n_0}$ pris au d\'ebut de la construction 
ne se rencontrent pas, la construction pr\'ec\'edente peut \^etre effectu\'ee s\'epar\`ement relativement aux deux strates 
 (voir figure \ref{thommather4}).  Cette construction est similaire \`a celle de M.-H. Schwartz \cite{MHS2}. 
 Ceci termine le cas de deux strates. 

\begin{figure}[h!]
\centering
\includegraphics[scale=0.7]{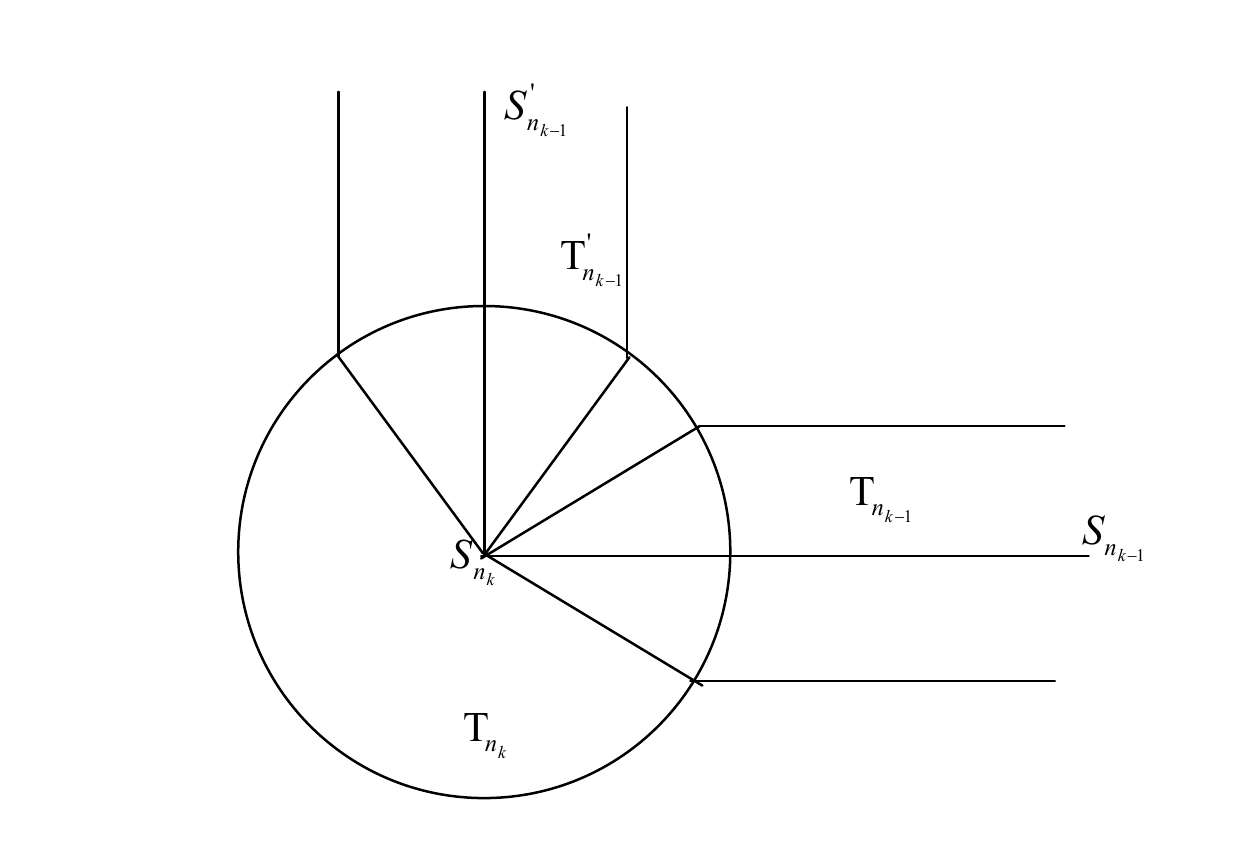}
\caption{Cas d'une strate dans le bord de deux strates.}
 \label{thommather4}
\end{figure}

Montrons maintenant que si la construction est faite pour une suite de strates 
$$S_{n_{k-1}}, \cdots, S_{n_1}, S_{n_0}$$
telles que 
$$S_{n_{k-1}} \subset  \cdots \subset{\overline{ S_{n_1}} } \subset {\overline{ S_{n_0}}},$$
alors, si $S_{n_k}$ est une strate telle que $S_{n_k} \subset {\overline{S_{n_{k-1}}}}$,
  la construction peut \^etre faite pour la suite de strates 
$$S_{n_k}, S_{n_{k-1}}, \cdots, S_{n_{1}}, S_{n_0}.$$
Consid\'erons une suite de fa\c cons 
$$\kappa_{n_k} \prec \kappa_{n_{k-1}} \prec \cdots \prec \kappa_{n_{1}} \prec \kappa_{n_0}$$
des strates correspondantes.
Par hypoth\`ese de r\'ecurrence, il existe des tubes ${\mathcal T}_{n_i}$ autour des strates 
$S_{n_i}$, ceci pour $i=0, \ldots, k-1$ ainsi que des fonctions $\pi_{n_i}$ et $\rho_{n_i}$
satisfaisant les conditions de Thom-Mather. 

Comme pr\'ec\'edemment (cf 1) du cas de deux strates), nous construisons d'abord  le voisinage tubulaire 
${\mathcal T}_{n_k}$ dans $V_{n_k} \bigcap_{i= k-1}^0 {\mathcal T}_{n_i}$. 
Pour un point $a \in \bigcap_{i= k-1}^0 {\mathcal T}_{n_i}$, 
nous  supposons donc construit un ``poly\`edre curviligne" $P_{k-1}$ de sommets le point $a$ et des points $a^{n_i}$
tels que $a^{n_i} \in S_{n_i}$, pour $i=k-1, \ldots, 0$.
Supposons le point $a$ suffisamment proche de $S_{n_k}$ ({\it i.e.} situ\'e dans le voisinage $V_{n_k}$)
et consid\'erons un point $a^{n_k} \in S_{n_k}$ proche de $a$. 

\begin{figure}[h!]
\centering
\includegraphics[scale=0.7]{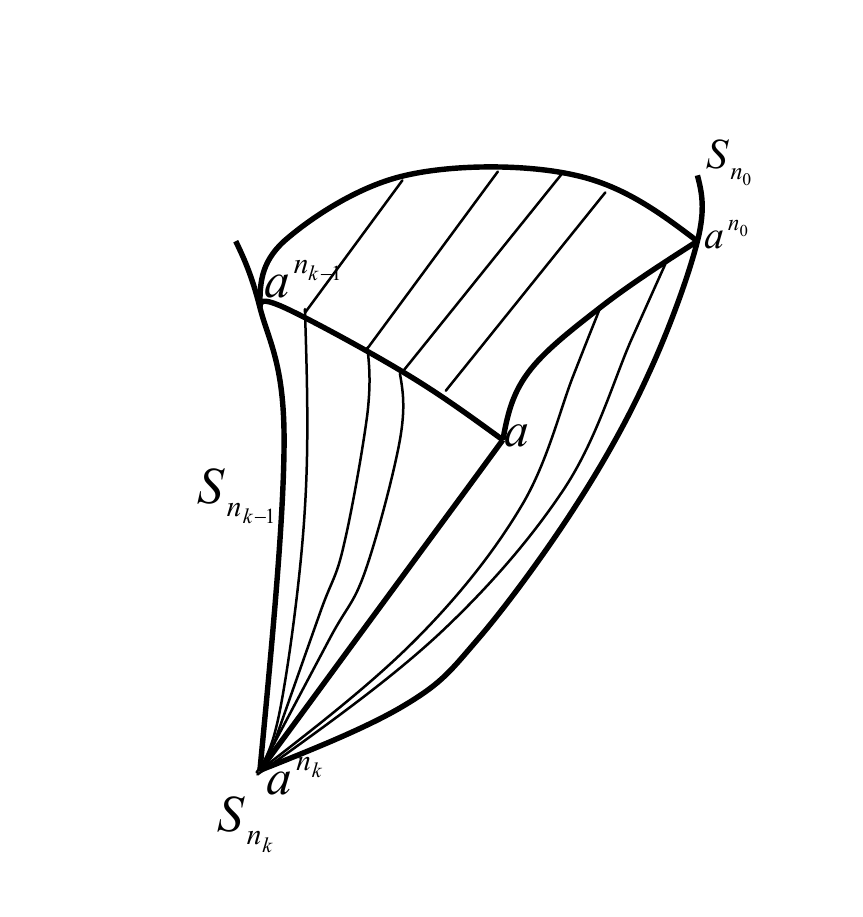}
\caption{La r\'ecurrence}
\label{thommather3}
\end{figure}

D'apr\`es le Lemme  \ref{lemmethuy3},
appliqu\'e \`a $a$ et $a^{n_k} $, il existe une courbe $\gamma^{n_k} (t_k)$ 
joignant $a$ et $a^{n_k} $, relativement \`a la fa\c con 
$\kappa_{n_k} $ (voir figure \ref{thommather3}). Comme pr\'ec\'edemment, chaque point $b$ d'une ar\^ete $a, a^{n_i}$ 
du poly\`edre curviligne peut \^etre joint au point
$a^{n_k} $ par une courbe du m\^eme type, qui est un rayon du voisinage tubulaire 
${\mathcal T}_{n_k}$ de $S_{n_k}$ que nous construisons. 
Lorsque le point $b$ tend vers $a^{n_i} \in S_{n_i}$, nous obtenons
une courbe $(a^{n_i} , a^{n_k})$ situ\'ee dans la strate $S_{n_i}$ et qui est \'egalement un rayon du voisinage tubulaire 
${\mathcal T}_{n_k}$ de $S_{n_k}$.  Les courbes obtenues d\'ependent toutes  du m\^eme param\`etre $t_k$. 
Nous obtenons ainsi un nouveau poly\`edre curviligne dont les sommets sont le point $a$ et
les points $a^{n_i}$ tels que $a^{n_i} \in S_{n_i}$, pour $i=k, \ldots, 0$.
Comme pr\'ec\'edemment, chaque point $a'$ de la courbe 
$(a, a^{n_k}) $ correspond \`a une valeur $t_k \in [1, +\infty)$ du param\`etre. 
En particulier, les points du poly\`edre $P_{k-1}$ correspondent \`a la valeur $t_k = 0$. Par d\'efinition, 
nous posons 
$$\rho_{n_k} (a) = \rho_{n_k} (a^{n_i}) \text{ pour } i= k-1, \ldots, 0$$
et tous les points du poly\`edre obtenu \`a partir des points de  $P_{k-1}$ et correspondant 
\`a la m\^eme valeur $t_k$ du param\`etre, ont m\^eme valeur de $\rho_{n_k}$. 

La construction de ${\mathcal T}_{n_k}$ dans $V_{n_k} \setminus \bigcap_{i= k-1}^0 {\mathcal T}_{n_i}$ 
ainsi que la d\'efinition des fonctions $\pi_{n_k}$ et $\rho_{n_k}$ est alors similaire au cas de deux 
strates effectu\'e ci-dessus. 

Il nous reste \`a pr\'eciser la construction des voisinages tubulaires dans les deux cas suivants :

1) cas o\`u il n'existe pas de relation entre les fa\c cons d'une strate et de la strate situ\'ee dans son bord.

2) cas o\`u l'on a plusieurs fa\c cons pour une m\^eme strate.

Dans le premier cas, puisque la stratification de $S_F$ d\'efinie par les fa{\c c}ons \'etoile est une stratification alg\'ebrique, alors il y a une relation entre des \'equations des strates adjacentes. En utilisant cette relation, la preuve est proc\'ed\'ee avec le m\^eme chemin. 

Dans le second cas, supposons qu'il existe $r$ fa\c cons pour une m\^eme strate $S_{\nu}$. 
Supposons qu'il existe deux suites de strates incluant $S_\nu$ pour lesquelles les fa\c cons $\kappa_\nu$ et $\kappa'_\nu$ 
sont diff\'erentes. Cela implique qu'il existe deux strates $S_{n_i}$ et $S'_{n_i}$ dans le bord de $S_\nu$ 
avec diff\'erentes fa\c cons. 
Les rayons de la partie du voisinage tubulaire 
${\mathcal T}_{\nu}$ situ\'ee dans ${\mathcal T}_{n_i}$ sont construits avec des courbes 
 $\gamma^{n_i}(u)$ tandis que les rayons de la partie du voisinage tubulaire 
${\mathcal T}_{\nu}$ situ\'ee dans ${\mathcal T}'_{n_i}$ sont construits avec des courbes 
 ${\gamma'}^{n_i}(u)$. 
Nous obtenons donc des syst\`emes de courbes diff\'erents pour chaque partie du voisinage tubulaire 
${\mathcal T}_{\nu}$ de $S_\nu$. Nous pouvons donc conclure comme dans le premier cas. En fait les courbes 
obtenues dans ce cas sont de la forme suivante :
$$\Phi(u, s) : = s_1 F(\gamma^{\nu}_{1}(u))  + \cdots + s_r F(\gamma^{\nu}_{r}(u)),$$
o\`u d'une part,  les  $r$ courbes $\gamma^{\nu}_{1}(u), \ldots, \gamma^{\nu}_{r}(u)$ tendant vers l'infini 
 et telles que leurs images tendent vers $a^{\nu}$ correspondent aux $r$ fa\c cons de $S_\nu$, et o\`u d'autre part, 
 $s_1  + \cdots + s_r = 1$ et $0 \le s_i \le 1$, pour $i = 1, \ldots, r$, seule l'une des valeurs $s_i$ valant $1$ dans 
 chaque ${\mathcal T}_{\nu} \cap {\mathcal T}_{n_i} $. 

\end{preuve}

\bibliographystyle{plain}

\end{document}